\documentstyle{article}

\input amssym.def
\input amssym

\let\got=\frak

\def\filename{newint.tex}

\def\hybrid{\topmargin 0pt      \oddsidemargin 0pt
        \headheight 0pt \headsep 0pt
        \textwidth 17.5cm
        \textheight 23cm
        \voffset=-0.5cm
        \hoffset=0.0cm
        \marginparwidth 0.0in
        \parskip 5pt plus 1pt   \jot = 1.5ex}
\catcode`\@=11
\def\marginnote#1{}

\newcount\hour
\newcount\minute
\newtoks\amorpm
\hour=\time\divide\hour by60
\minute=\time{\multiply\hour by60 \global\advance\minute by-\hour}
\edef\standardtime{{\ifnum\hour<12 \global\amorpm={am}%
        \else\global\amorpm={pm}\advance\hour by-12 \fi
        \ifnum\hour=0 \hour=12 \fi
        \number\hour:\ifnum\minute<10 0\fi\number\minute\the\amorpm}}
\edef\militarytime{\number\hour:\ifnum\minute<10 0\fi\number\minute}

\def\draftlabel#1{{\@bsphack\if@filesw {\let\thepage\relax
   \xdef\@gtempa{\write\@auxout{\string
      \newlabel{#1}{{\@currentlabel}{\thepage}}}}}\@gtempa
   \if@nobreak \ifvmode\nobreak\fi\fi\fi\@esphack}
        \gdef\@eqnlabel{#1}}
\def\@eqnlabel{}
\def\@vacuum{}
\def\draftmarginnote#1{\marginpar{\raggedright\scriptsize\tt#1}}

\def\draft{\oddsidemargin -0.1truein
        \def\@oddfoot{\sl preliminary draft\ {\tt\filename} \hfil
        \rm\thepage\hfil\sl\today\quad\militarytime}
        \let\@evenfoot\@oddfoot \overfullrule 3pt
        \let\label=\draftlabel
        \let\marginnote=\draftmarginnote
   \def\@eqnnum{{\rm (\theequation)}\rlap{\kern\marginparsep\tt\@eqnlabel}%
\global\let\@eqnlabel\@vacuum}  }

\newcommand{\CC}{{\Bbb C}}
\newcommand{\NN}{{\Bbb N}}

\newcommand{\ZZ}{{\Bbb Z}}
\let\z=\ZZ

\font\fiveeufm=cmmib5  

\font\yao=cmmi12 
\def\bep{{\textfont1=\yao\scriptfont1=\teneuf
\scriptscriptfont1=\fiveeufm
\hbox{$\mathsurround=0pt\displaystyle\epsilon$}}
}

\newdimen\linethick  \linethick=0.4pt
\newdimen\hboxitspace    \hboxitspace=5pt
\newdimen\vboxitspace    \vboxitspace=5pt

\def\fr#1{%
\beq\new
\vcenter{
\hrule height\linethick
           \hbox{\vrule width\linethick
                 \kern\hboxitspace
                 \vbox{\kern\vboxitspace
                       \hbox{$\begin{array}{c}\displaystyle#1
          \end{array}$}%
                       \kern\vboxitspace}%
                 \kern\hboxitspace
                 \vrule width\linethick}%
           \hrule height\linethick}%
\eeq}

\newdimen\Squaresize \Squaresize=14pt
\newdimen\Thickness \Thickness=0.5pt

\def\Square#1{\hbox{\vrule width \Thickness
   \vbox to \Squaresize{\hrule height \Thickness\vss
      \hbox to \Squaresize{\hss#1\hss}
   \vss\hrule height\Thickness}
\unskip\vrule width \Thickness}
\kern-\Thickness}

\def\Vsquare#1{\vbox{\Square{$#1$}}\kern-\Thickness}

\def\numberbysection{\@addtoreset{equation}{section}
        \def\theequation{\thesection.\arabic{equation}}}
\numberbysection

\renewcommand{\theequation}{\thesection.\arabic{equation}}
\newcommand{\l@qq}[2]{\addvspace{2em}
 \hbox to\textwidth{\hspace{1em}\bf #1 \dotfill #2}}

\def\appname{Appendix}
\newcounter{app}
\def\theapp{\Alph{app}}
\def\app{\par
   \addvspace{4ex}
   \@afterindentfalse
  \secdef\@app\@dapp}
\def\@app[#1]#2{\ifnum \c@secnumdepth >\m@ne
        \refstepcounter{app}
        \addcontentsline{toc}{app}{\theapp
        \hspace{1em}#1}\else
      \addcontentsline{toc}{app}{ #1}\fi
   {\parindent \z@ \raggedright
    \Large \bf \appname~\theapp .
   \Large  \bf \hspace{1em}    #2}\nobreak
   \vskip -2ex   \noindent
\setcounter{equation}{0}
\def\theequation{\Alph{app}.\arabic{equation}}}
\def\@dapp#1{%
{\parindent \z@ \raggedright  \bf #1}\par\nobreak}
\def\l@app#1#2{\addpenalty{\@secpenalty}%
   \addvspace{1em plus\p@}%
   \begingroup
   \@tempdima 3em
     \parindent \z@ \rightskip \@pnumwidth
     \parfillskip -\@pnumwidth
     { \bf
     \leavevmode
     #1\hfil \hbox to\@pnumwidth{\hss #2}}\par
     \nobreak
   \endgroup}
\newcounter{sapp}[app]
\def\thesapp{\Alph{app}.\arabic{sapp}}
\def\sapp{\par
   \addvspace{4ex}
   \@afterindentfalse
  \secdef\@sapp\@dsapp}
\def\@sapp[#1]#2{\ifnum \c@secnumdepth >\m@ne
        \refstepcounter{sapp}
        \addcontentsline{toc}{sapp}{\thesapp
        \hspace{1em}#1}\else
      \addcontentsline{toc}{sapp}{ #1}\fi
   {\parindent \z@ \raggedright
    \large \bf \thesapp
   \large  \bf \hspace{1em}    #2}\nobreak
   \vskip 4ex   \noindent
\def\theequation{\Alph{app}.\arabic{equation}}}
\def\@dsapp#1{%
{\parindent \z@ \raggedright  \bf #1}\par\nobreak}
\def\l@sapp#1#2{\addpenalty{\@secpenalty}%
   \begingroup
   \@tempdima 3em
     \parindent \z@ \rightskip \@pnumwidth
     \parfillskip -\@pnumwidth
     { \hspace{1em}
     \leavevmode
     #1 \hfil \dotfill \hbox to\@pnumwidth{\hss #2}}\par \nobreak
     \endgroup}

\def\titlepage{\@restonecolfalse\if@twocolumn\@restonecoltrue\onecolumn
     \else \newpage \fi \thispagestyle{empty}\c@page\z@
        \def\thefootnote{\fnsymbol{footnote}} }

\def\endtitlepage{\if@restonecol\twocolumn \else  \fi
        \def\thefootnote{\arabic{footnote}}
        \setcounter{footnote}{0}}  
\relax

\hybrid

\newtoks\@stequation

\def\subequations{\refstepcounter{equation}%
  \edef\@savedequation{\the\c@equation}%
  \@stequation=\expandafter{\theequation}
  \edef\@savedtheequation{\the\@stequation}
  \edef\oldtheequation{\theequation}%
  \setcounter{equation}{0}%
  \def\theequation{\oldtheequation\alph{equation}}}

\def\endsubequations{%
  \setcounter{equation}{\@savedequation}%
  \@stequation=\expandafter{\@savedtheequation}%
  \edef\theequation{\the\@stequation}%
  \global\@ignoretrue}

\makeatletter
\newdimen\normalarrayskip              
\newdimen\minarrayskip                 
\normalarrayskip\baselineskip
\minarrayskip\jot
\newif\ifold             \oldtrue            \def\new{\oldfalse}
\def\arraymode{\ifold\relax\else\displaystyle\fi} 
\def\eqnumphantom{\phantom{(\theequation)}}     
\def\@arrayskip{\ifold\baselineskip\z@\lineskip\z@
     \else
     \baselineskip\minarrayskip\lineskip1\baselineskip\fi}

\def\@arrayclassz{\ifcase \@lastchclass \@acolampacol \or
\@ampacol \or \or \or \@addamp \or
   \@acolampacol \or \@firstampfalse \@acol \fi
\edef\@preamble{\@preamble
  \ifcase \@chnum
     \hfil$\relax\arraymode\@sharp$\hfil
     \or $\relax\arraymode\@sharp$\hfil
     \or \hfil$\relax\arraymode\@sharp$\fi}}

\def\@array[#1]#2{\setbox\@arstrutbox=\hbox{\vrule
     height\arraystretch \ht\strutbox
     depth\arraystretch \dp\strutbox
     width\z@}\@mkpream{#2}\edef\@preamble{\halign \noexpand\@halignto
\bgroup \tabskip\z@ \@arstrut \@preamble \tabskip\z@ \cr}%
\let\@startpbox\@@startpbox \let\@endpbox\@@endpbox
  \if #1t\vtop \else \if#1b\vbox \else \vcenter \fi\fi
  \bgroup \let\par\relax
  \let\@sharp##\let\protect\relax
  \@arrayskip\@preamble}
\def\eqnarray{\stepcounter{equation}%
              \let\@currentlabel=\theequation
              \global\@eqnswtrue
              \global\@eqcnt\z@
              \tabskip\@centering                      
              \let\\=\@eqncr
              $$%
            \halign to \displaywidth  \bgroup
             \eqnumphantom \@eqnsel
      \hskip\@centering                               
    $\displaystyle  \tabskip\z@ {##}$%
    &\global\@eqcnt\@ne \hskip 2\arraycolsep
         $ \displaystyle  \arraymode{##}$\hfil
    &\global\@eqcnt\tw@ \hskip 2\arraycolsep
         $\displaystyle\tabskip\z@{##}$\hfil
         \tabskip\@centering
    &{##}\tabskip\z@\cr}
\makeatother
\newtheorem{th}{Theorem}
\newtheorem{prop}{Proposition}[section]           

\newtheorem{lem}{Lemma}[section]

\def\bea{\begin{eqnarray}}
\def\eea{\end{eqnarray}}

\def\beq{\begin{equation}}
\def\eeq{\end{equation}}
\def\be{\beq\new\begin{array}{c}}
\def\ee{\end{array}\eeq}
\def\bse{\begin{subequations}}                
\def\ese{\end{subequations}}                 %
\def\vbn{$\raisebox{0.1pt}{$\,\hspace{2pt}\llcorner\!$}\leaders
\hrule\hfill\leaders\hrule\hfill\raisebox{0.1pt}{$\!\lrcorner$}\;\,$}
\def\con#1{\mathsurround=0pt
\mathop{\vtop{\ialign{##\crcr$\hfil\displaystyle{#1}\hfil$\crcr
\noalign{\kern1pt\nointerlineskip}\vbn\crcr\noalign{\kern1pt}}}}\limits}
\def\stackreb#1#2{\mathrel{\mathop{#2}\limits_{#1}}}

\def\tr{\triangleright}                                  
\def\tl{\triangleleft}

\def\jo{\mathrel{\mkern-4mu}}
\def\sem{\mathsurround=0pt
\mathrel{\raise1.4pt\hbox{$\scriptscriptstyle>$}}\jo\mathrel\tl}
\def\mes{\mathsurround=0pt
{\mathrel\tr\jo\mathrel{\raise1.4pt\hbox{$\scriptscriptstyle <$}}}}

\def\]{]\raise-2pt\hbox{$_\ast$}}

\def\op#1{\raise-6pt\hbox{$\stackrel{\displaystyle\oplus }{\scriptstyle
#1}$}\;}
\def\oop#1#2{\raise-6pt\hbox{$\stackrel{#2}{\stackrel{\displaystyle\oplus
}{\scriptstyle #1}}$\;}}

\def\al{\alpha}
\def\a{\alpha}
\def\b{\beta}
\def\la{\lambda}

\def\xp{e}
\def\xm{f}
\def\co{\mathsurround=0pt\e\raise-4pt
\hbox{${\!{_{_\cf}\!{\scriptscriptstyle\cal H}}}$}}
\def\dco{\mathsurround=0pt\bep\raise-4pt
\hbox{${\!\!{_{_\cf}\!\!{\scriptscriptstyle\cal H^{\!\ast}}}}$}}
\def\<{\langle}
\def\>{\rangle}
\def\ov{\overline}

\def\g{{\bf g}}
\def\ggg{\widehat{{\bf g}}}
\def\hh{{\bf h}}
\def\hhh{\widehat{\hh}}

\def\cu{{\cal U}}
\def\kr{{\cal R}}

\def\R{\ov\kr}

\def\cf{{\cal F}}
\def\1{1\!\!1}
\def\stack#1#2{\raise0.7pt\hbox{$\mathrel{\mathop{#2}\limits^{#1}}$}}
\def\m{\raise5pt\hbox{$\scriptstyle -1$}}
\def\st#1#2{\raise1pt\hbox{$\stackrel{#1}{#2}$}}

\def\fd{\cu\!\mathsurround=0pt
\raise-7pt\hbox{$\stackrel{\displaystyle\otimes}
{\scriptscriptstyle \cal R^{\!\scriptscriptstyle(\!-\!)}}$}\cu}

\def\r#1{(\ref{#1})}

\def\ot{\otimes}
\def\a{\alpha}  \def\d{\delta}  \def\la{\lambda}
\def\b{\beta} \def\ga{\gamma} 

\def\sldva{{\frak sl}_2}
\def\D{\Delta}

\newcommand{\Uqg}{{U_{q}^{}(\widehat{{\bf g}})}}
\newcommand{\UqgD}{{U_{q}^{(D)}(\widehat{{\bf g}})}}
\newcommand{\Uqdva}{{U_{q}^{}(\widehat{{\bf sl}}_{2})}}
\newcommand{\UqdvaD}{{U_{q}^{(D)}(\widehat{{\bf sl}}_{2})}}

\def\Pqexp{{\cal P}{\overrightarrow{\exp}}_{\{D_{I}\}}}
\def\K{{\cal K}}

\let\bn=\be
\let\ed=\ee
\let\rf=\r
\def\sk#1{\left(#1\right)}

\def\I#1{I^{(#1)}}
\def\e#1{e^{({#1})}}
\def\eee#1{\tilde{e}^{({#1})}}
\def\f#1{f^{({#1})}}
\def\T#1{t^{({#1})}}

\def\krr#1{\ov\kr^{(#1)}}
\def\dz{\underline{dz}}
\def\dzz#1{\underline{dz_{#1}}}
\def\ddz#1{\frac{dz_{#1}}{z_{#1}}}

\def\res#1{\stackreb{#1}{\mbox{\rm res}}}
\def\lmt#1{\stackreb{#1}{\mbox{\rm lim}}}
\let\Res=\res

\def\krr#1{\ov{{\cal R}}^{(#1)}}

\def\II{\check{I}}
\def\io{\iota}
\def\DD{{D}}
\begin{document}
\thispagestyle{empty}
\setcounter{page}1

\begin{center}
\hfill ITEP-TH-67/99\\
\hfill math.QA/0008226\\
\bigskip\bigskip
{\Large\bf Integral Presentations for the Universal ${\cal R}$-matrix}\\
\bigskip
\bigskip
{\bf
J. Ding$^{*}$\footnote{E-mail: Jintai.Ding@math.uc.edu},\
S. Khoroshkin$^{\star}$\footnote{E-mail: khor@heron.itep.ru},\
S. Pakuliak$^{\star\bullet\circ}$}\footnote{E-mail: pakuliak@thsun1.jinr.ru}\\
\bigskip
$^{*}$
{\it Department of  Mathematical Sciences, University of Cincinnati,\\
 PO Box 210025, Cincinnati, OH 45221-0025,  USA}\\
$^\star$
{\it Institute of Theoretical \& Experimental Physics, 117259 Moscow,
Russia}\\
$^\bullet$
{\it Bogoliubov Laboratory of Theoretical Physics, JINR,
141980 Dubna, Moscow region, Russia}\\
$^\circ$
{\it Bogoliubov Institute for Theoretical Physics, 252143 Kiev, Ukraine}
\bigskip
\bigskip
\end{center}
\medskip
\hfill{{\sl to the memory of Denis Uglov \hspace{2cm}}}
\bigskip
\begin{abstract}
We present an integral formula for the universal ${\cal R}$-matrix of
 quantum
affine algebra $\Uqg$ with 'Drinfeld comultiplication'. We show
that the properties of the universal ${\cal R}$-matrix
 follow from the factorization
properties of the cycles in proper configuration spaces. For general
${\bf g}$ we
 conjecture that such cycles exist and unique. For $\Uqdva$ we
 describe precisely the cycles and present a new simple expression for
 the universal ${\cal R}$-matrix as a result of calculation of corresponding
integrals.
  \end{abstract} 

\footnotesize
\tableofcontents
\newpage
\normalsize

\section{Introduction}
The Yang-Baxter (YB) equation is one of cornerstones in the investigations
 of quantum integrable systems. The most important solutions to the YB equation
 were found by Yang, Baxter \cite{Baxter} and Zamolodchikovs
\cite{Zamolodchikov}. The mathematical background to the application of the
YB equation was established in the theory of quantum groups by Drinfeld
\cite{D} and Jimbo \cite{J}, based on the quantum inverse scattering
method developed by Leningrad school \cite{Faddeev}.

The main object in this theory is a quantum group which is a Hopf algebra
deformation of the universal enveloping algebra  of contragredient Lie algebra.
 This Hopf algebra comes together with the universal $R$-matrix which acts
 in a tensor category of representations of quantum group and any such a
 specialization of the universal $R$-matrix provides a solution of the YB
equation. The deformations of affine Lie algebras : Yangians, quantum affine
algebras and their elliptic analogs are especially important: they produce
$R$-matrices depending on a spectral parameter, in particular the main
solutions of the
 YB equation listed above.

However, $R$-matrices and fundamental $L$-operators appear in quite
different manner in physical and mathematical literature. In the works, written
by physicists, we have their expressions written as ordered exponential integrals
 like a series of iterated integrals \cite{BLZ}, while in the mathematical
literature all the known expressions can be reduced to products of
certain $q$-exponential factors \cite{KT1}.  This reflects the structure
 of quantum affine algebras: they deform the structure of contragredient
 Lie algebra for affine Lie algebras and do not touch actually their
functional realization.

Drinfeld suggested another ('new') realization of quantum affine algebras
 and of Yangians \cite{D1}, which, together
with naturally attached coproduct
 structure, can be regarded as a Hopf algebra deformation of affine Lie
 algebras in a sense that it deforms their functional realization.
This approach was successfully exploited for the study of finite-dimensional
representations of quantum affine algebras and of Yangians (see \cite{Chari}
 and references therein) and for bosonization of their infinite-dimensional
 representations
 \cite{FJ}.

In \cite{DK} a functional version of braid group technique was developed for
quantum affine algebras and for the doubles of the Yangians. In particular,
the currents corresponding to nonsimple roots of underlying Lie algebras were
 defined and a presentation of the universal $R$-matrix which uses
 multiple contour integrals was given.

In this note we develop further the ideas of \cite{DK}.
The first main result is a geometric reformulation of the
properties of the universal $R$-matrix for quantum affine algebras equipped
with so called Drinfeld comultiplication. Namely, we attach to any simple
Lie algebra $\g$ a system of configuration spaces, which are complex linear
spaces without shifted diagonals and coordinate hyperplanes. We claim that
the defining properties of the universal $R$-matrix can be reduced
to certain factorization properties
 for the homology classes of the cycles in configuration space.
The cycles are actually relative: we are interested in their pairing
with the forms satisfying certain zero residues conditions. In algebra, these
vanishing conditions
 come from the Serre relations.
 Thus we have a purely homological problem
 which resembles the constructions of factorizable sheaves
\cite{BFS} and  of factorizable $D$-modules \cite{KS}.
In general case we do not give here a construction of factorizable system
 of cycles, just formulate a conjecture of its existence.

Second result is the explicit  calculation of the $R$-matrix for $\Uqdva$.
 We suggest in this case the factorizable system of cycles in configuration
space, generalizing the ideas of \cite{DK} where it was done for $|q|>1$
 and calculate the corresponding integrals. The calculations
 go as follows: we take some residues to derive the recurrence relations
 for $n$-th fold integral, which lead to a simple linear differential equation
 for the $R$-matrix. Its solution has a form of vertex operator over
 total integrals of residue currents. These currents appeared already in the
study of integrable representations  of $\Uqdva$ \cite{DM}: the annihilator of
level $k$ integrable representations is described in their terms. As a
consequence, the $R$-matrix in this case
contains only finitely many currents which
contribute
 into vertex operator, for level one it is just an exponent of total integral
 of the first one. An application of these results to quantum affine algebra
$\Uqdva$ with traditional comultiplication is given in \cite{DKP}.

\section{Quantized current algebras}

 By quantized current algebra we mean a completion of quantum affine algebra
  in Drinfeld ''new'' realization equipped with coproduct structure,
 introduced
 by Drinfeld in the deriving of this realization \cite{D1}. The completion is
done   with respect to a minimal topology, in which the action in the graded
modules with bounded from above degree is continuous.
\subsection{''New'' realization of quantum affine algebras}

Let $\g$ be a simple Lie algebra. Following Drinfeld \cite{D1}, we can
 describe
quantum affine algebra $\Uqg$ as an algebra generated by the central
element $c$, gradation operator $d$ and by the modes of the currents
 $e_{\pm\al_i}(z)$ and $\psi^\pm_{\al_i}(z)$,
\be
e_{\pm\al_i}(z)=\sum_{k\in\z}e_{\pm\al_i,k}z^{-k},\qquad
\psi^\pm_{\a_i}(z)=k_{\al_i}^{\pm 1}\exp\left(\pm(q-q^{-1})\sum_{n>0}
a_{i,\pm n}z^{\mp n}\right),
\label{Dr-g}
\ee
where
$\al_i$, $i=1,\ldots,r$ are  positive simple roots of the algebra $\g$ of the
rank $r$,
and $(\al_i,\al_j)$ is symmetrized Cartan
matrix.

 The generating functions of the quantum current algebra
$\Uqg$  satisfy the following defining relations:
\be\label{d00}
q^d x(z) q^{-d}=x(qz),\quad\mbox{for}\quad
x=e_{\pm\al_i},\ \psi^\pm_{\al_i}\ ,
\ee
\bn
(z-q^{\pm({\al_i},{\al_j})}w)e_{\pm\al_i}(z)e_{\pm\al_j}(w)=
 e_{\pm\al_j}(w)e_{\pm\al_i}(z)(q^{\pm({\al_i},{\al_j})}z-w)\ ,
\label{1}
\ed
\bn
\frac{(q^{\pm c/2}z-q^{\pm({\al_i},{\al_j})}w)}
{(q^{\pm({\al_i},{\al_j})\pm c/2}z-w)}
 \psi_{\al_i}^+(z)e_{\pm\al_j}(w)=
e_{\pm\al_j}(w)\psi_{\al_i}^+(z)\ ,
\label{3}
\ed
\bn
\psi_{\al_i}^-(z)e_{\pm\al_j}(w)=\frac{(q^{\pm({\al_i},{\al_j})\mp c/2}z-w)}
{(q^{\mp c/2}z- q^{\pm({\al_i},{\al_j})}w)}
e_{\pm\al_j}(w)\psi_{\al_i}^-(z)\ ,
\label{4}
\ed
\bn
\frac{(z-q^{({\al_i},{\al_j})-c}w)(z-q^{-({\al_i},{\al_j})+c}w)}
{(q^{({\al_i},{\al_j})+c}z-w)(q^{-({\al_i},{\al_j})-c}z-w)}
\psi_{\al_i}^+(z)\psi_{\al_j}^-(w)=\psi_{\al_j}^-(w)\psi_{\al_i}^+(z)\ ,
\label{7}
\ed
\bn
\psi_{\al_i}^\pm(z)\psi_{\al_j}^\pm(w)=
 \psi_{\al_j}^\pm(w)\psi_{\al_i}^\pm(z)\ ,
\label{7a}
\ed
\bn
[e_{\al_i}(z),e_{-\al_j}(w)]=
 \frac{\delta_{{\al_i},{\al_j}}}{q_{}-q^{-1}_{}}\left( \delta(z/q^cw)
\psi^+_{\al_i}(zq^{-c/2})-\delta(zq^c/w)\psi^-_{\al_i}(wq^{-c/2})\right)
\label{10}
\ed
and the Serre relations:
\bn
\sum_{r=0}^{n_{ij}}(-1)^r\left[{n_{ij}\atop r}\right]_{q_{\a_i}}
{\rm Sym}_{z_1,...,z_{n_{ij}}}e_{\pm\al_i}(z_1)\cdots e_{\pm\al_i}(z_r)
e_{\pm\al_j}(w)e_{\pm\al_i}(z_{r+1})\cdots e_{\pm\al_i}(z_{n_{ij}})=0
\label{serre1}
\ed
 for any simple roots ${\al_i}\not={\al_j}$.
Here $\delta(z)=\sum_{k\in \z}z^k$, $n_{ij}=1-a_{ij}$, where $a_{ij}$ is
$(i,j)$ entry of Cartan matrix of $\g$,
$\left[{n\atop k}\right]_q=\frac{[n]_q!}{[k]_q![n-k]_q!}$,
$[n]_q!=[1]_q[2]_q\cdots [n]_q$, $[n]_q=\frac{q^n-q^{-n}}{q-q^{-1}}$,
 $q_\a=q^{\frac{(\a,\a)}{2}}$.

 The formulas \r{3}--\r{7} can be rewritten also using the generators
$a_{i,n}$ introduced by \r{Dr-g}:
\be\label{34a}
[a_{i,n},e_{\pm{\al}_j}(z)]
 =\pm\frac{[({\al_i},{{\al_j}})n]_q}{n}q^{\mp c|n|/2}z^n e_{\pm{\al_j}}(z)\ ,
\ee
\be\label{777}
[a_{i,n},a_{j,m}]
 =\d_{n+m,0}\frac{[({\al_i},{{\al_j}})n]_q[cn]_q}{n}\ .
\ee

We assign to this algebra two Drinfeld type Hopf structures:
\bse
\label{Pa}
\bn
\Delta^{(1)} e_{{\a_i}}(z)=
e_{{\a_i}}(z)\ot 1 +\psi_{{\a_i}}^-(zq^{c_1/2})\ot e_{{\a_i}}(zq^{c_1})\ ,
\label{Pa1}
\ed
\bn
\Delta^{(1)} e_{-{\a_i}}(z)=
1\ot e_{-{\a_i}}(z) +e_{-{\a_i}}(zq^{c_2})\ot \psi_{{\a_i}}^+(zq^{c_2/2})\ ,
\label{Pa2}
\ed
\ese
\bse\label{Paa}
\bn
\Delta^{(2)} e_{{\a_i}}(z)=
e_{{\a_i}}(z)\ot 1 +\psi_{{\a_i}}^+(zq^{-c_1/2})\ot e_{{\a_i}}(zq^{-c_1})\ ,
\label{Paa1}
\ed
\bn
\Delta^{(2)} e_{-{\a_i}}(z)=
1\ot e_{-{\a_i}}(z) +e_{-{\a_i}}(zq^{-c_2})\ot \psi_{{\a_i}}^-(zq^{-c_2/2})\ ,
\label{Paa2}
\ed
\ese
\bn
\Delta^{(1)} \psi_{{\a_i}}^\pm(z)=
\Delta^{(2)} \psi_{{\a_i}}^\pm(z)=\psi_{{\a_i}}^\pm(zq^{\pm c_2/2})\ot
\psi_{{\a_i}}^\pm(zq^{\mp c_1/2})\ ,
\label{Paa3}
\ed
with corresponding expressions for counit and antipode maps.
One can see that the comultiplication operators map an algebra $\Uqg$
  to a completed tensor product $\Uqg\hat{\otimes}\Uqg$. One
possibility to describe it is to  use in the
 right hand side the
Taylor extension of $\Uqg\otimes \Uqg$, defined in \cite{KT}.
Below we define instead a completion in a weaker topology, adapted for the
 action on highest degree representations.

\subsection{The completion}

 Let $U=\Uqg$.
For any $x\in U$,  we denote
$|x|<k$ if $x$ can be presented as a noncommutative polynomial of degree
less then $ k$ over the generators of the algebra, which are listed above.

The grading element $d$ defines a
 gradation on $U$: $\deg x=r$, if $[d,x]=rx$.
 Let $U_r$ be a linear subspace of $U$, generated by the elements of
degree $r$ and  $U^k$ be a linear subspace of $U$ generated by all $x\in U$,
such that $|x|+|\deg x|<k$. It is clear, that $\CC 1=U^1\subset U^2\subset
 U^3\subset...$ and $\cup U_k=U$.

Define  $U_r^k$ as
$$
U_r^k=U^k\cap\sum_{s\geq r}U\cdot U_s.
$$
We denote by $U^{<}$ a topological vector space $U$ with $U_r^k$
 being basic open neighborhoods of zero, and by $\ov{U}^<$
 the corresponding completion.

  One can note, that this topology
 can be defined in equivalent way by the system
of open sets
$${}_rU^k=U^k\cap\sum_{s\leq r}U_s\cdot U_.$$

 Let $W=U\ot U$. It is bigraded algebra, where $\deg_1$ is defined by the
 action
  of $d\ot 1$, and $\deg_2$ is defined by the
 action of $1\ot d$. We consider
 $x_k\ot 1$ and $1\ot x_k$ as generators of $W$, where $x_k$ are
generators of
 $U$. In this setting we define $W^k$ as linear span of all $x\in W$
such that
  $|\deg_1x|+|\deg_2x|+|x|<k$ and
 \bn
  W_r^k=W^k\cap \left(U_r^k\otimes U+U\ot U_r^k\right).
  \label{wrk}
  \ed
 We denote by $W^<$ the corresponding topological vector space
  and by $\ov{W}^<$ its completion.
 Let $\Delta^{(1)}$ and $\Delta^{(2)}$ \r{Pa} and \r{Paa}
 be two types of Drinfeld
 comultiplications in $U_q(\widehat{\frak{g}})$.
  \begin{prop}
(\i) The multiplication
 $$m: W^<\longrightarrow U^{<}$$
  is continuous map;

\noindent (\i\i)  The comultiplications
 $$ \Delta^{(1)},\ \Delta^{(2)}:U^{<}\longrightarrow
  \ov{W}^< $$
 are continuous maps;

\noindent (\i\i\i) The completed algebra $\ov{U}^<$ acts on the highest
degree modules over  $U$;

\noindent (\i v) The completed algebra $\ov{W}^<$ acts on the tensor
 products of highest degree modules over  $U$.
\end{prop}

Here by the highest degree representation we mean graded (over operator $d$)
 representation
$V$, such that the graded components of $V$ vanish for big enough degree:

$$V=\sum_{k<N}V_k.$$
For instance, all highest weight representations of $\Uqg$ are of highest
degree.

Actually, `$<$' is the smallest topology, in which the action
on any highest degree representation is continuous. We will use this
  property further checking the assertions on highest degree modules.
 Moreover,
  as we see from the proposition, the completion $\ov{U}^<$ is well defined
  topological Hopf algebra in the category of highest degree
 representations.

In quite analogous manner, we can define the
topological Hopf algebra
$U^>$ by means of open sets
\bn
  \tilde{U}_r^k=U^k\cap\sum_{s\leq r}U\cdot U_s.
  \label{turk}
  \ed
 Then its completion   $\ov{U}^>$ is well defined
 topological Hopf algebra in the category of lowest degree representations.
Denote further the completed  algebra $\ov{U}^<$ equipped with
comultiplication $\Delta^{(1)}$ as $\UqgD$.

The described completion $\UqgD$ of the algebra $\Uqg$
changes its algebraical properties.
In particular, it cannot be any more treated as a double of its Hopf
subalgebra, generated by $e_{-\a_i,n}$, $n\in\ZZ$ and $a_{i,k}$, $k\geq
0$, since in the dual subalgebra the topology should be reversed. In this
point we change the ideology of \cite{DK}. However, unlike the known
examples, the Hopf algebra $\UqgD$
   possesses the universal $R$-matrix, which acts in tensor
category of its representations.

\section{Universal $\kr$-matrices}

 The theory of Cartan-Weyl basis for quantum affine algebras, developed
 in \cite{TK}, \cite{KT}, \cite{KT1},  \cite{Be} gives
 possibility to write down the universal $R$-matrix for quantum affine
 algebras equipped with Drinfeld's comultiplication, in a form of infinite
 product over Cartan-Weyl generators. It was done in \cite{KT1}. This
  formula is quite effective for $\Uqdva$ but does not look nice for
 higher rank. We remind it below. Next, we analyze the integral form of the
  $R$-matrix for the completed algebra $\UqdvaD$, presented in \cite{DK}
and observe the geometric properties of
 the cycles  in certain configuration space, responsible for the properties
  of the universal $R$-matrix. As a corollary, we suggest the generalization
  of the formula in \cite{DK} for other $\UqgD$. The next section is devoted to
 explicit calculations in $\UqdvaD$.

 \subsection{The multiplicative formula \cite{KT1}}

 In the study of twistings of quantum affine algebra by means of the
elements  of the braid group the multiplicative formula for
the universal $R$-matrix for Drinfeld comultiplications
was suggested in \cite{KT1}. It was written by means of Cartan-Weyl
generators and can be described as follows.

 Let $e_{\pm\a_i},
k_{\a_i}^{\pm 1}$, $i=0,1,...,r$ be Chevalley generators of $\Uqg$, where
$r= {\rm rank}\ \g$. For their connection to current generators
$e_{\pm\a_i,n}$,
 $a_{i,n}$ see \cite{D1}. Let $\underline{\Delta}$ be a system
of positive roots of affine Lie algebra $\ggg$. Let $\a_0$ be an affine root
 and
$\delta$ be a minimal imaginary root.
 By a normal ordering of the root
system $\underline{\Delta}$ we mean a total linear order $<$ of all the roots
satisfying the condition: for any two roots $\a$, $\b$, such that at least
 one of them is real and $\a+\b$ is a root only two possibility of mutual
 position of
 $\a$, $\b$, $\a+\b$ can occur: $\a<\a+\b<\b$ or $\b<\a+\b<\a$.
Let $<$ be arbitrary normal ordering
of the root system $\underline{\Delta}$, satisfying additional restriction
 $\a_k<\delta<\a_0$ for any $k=1,...,r$.
The Cartan-Weyl generators are constructed as successive $q$-commutators by
induction on the normal order. An induction step looks as follows. We put
 \bn
e_{\a+\b}=[e_\a,e_\b]_q,\qquad e_{-\a-\b}=[e_{-\b},e_{-\a}]_{q^{-1}}
\label{induction}
\ed
if $e_{\a}$, $e_{-\a}$, $e_{\b}$, $e_{-\b}$ are already defined,
$\a<\a+\b<\b$ and $[\a,\b]$ is a minimal segment, including $\a+\b$, that is
if there are other roots $\a'$, $\b'$, such that $\a'+\b'=\a+\b$ and
$\a'<\a+\b<\b'$ then either $\a'<\a$ or $\b<\b'$. The $q$-commutators
$[e_\a,e_\b]_q$ and $[e_{-\b},e_{-\a}]_{q^{-1}}$ mean
$$[e_\a,e_\b]_q=e_\a e_\b-q^{(\a,\b)}e_\b e_\a,\qquad
[e_{-\b},e_{-\a}]_{q^{-1}}= e_{-\b}e_{-\a}-q^{-(\a,\b)}e_{-\a}e_{-\b}.$$

More precisely, we start from the simple roots and use the induction
procedure \r{induction} for the construction of new {\it real} roots until
we get the root vectors $e_{\pm(\delta-\a_i)}$, $i=1,...,r$. Then we put
$$
{e_{\delta}^{(i)}}=[e_{{\alpha}_{i}}, e_{\delta -{\alpha}_{i}}]_{q}
,\qquad
{e_{-\delta}^{(i)}}=[ e_{\delta -{\alpha}_{i}},e_{{\alpha}_{i}}]_{q^{-1}}  ,
$$
$$
e_{n\delta + \alpha_{i}}=(-1)^n([(\alpha_i,\alpha_i)]_q)^{-n}
(\tilde{\rm ad} \; e_{\delta}^{(i)})^{n}e_{{\alpha}_{i}} ,
\qquad
e_{-n\delta - \alpha_{i}}=([(\alpha_i,\alpha_i)]_q)^{-n}
(\tilde{\rm ad} \; e_{-\delta}^{(i)})^{n}e_{-{\alpha}_{i}}  ,
$$
$$
e_{(n+1)\delta -\alpha_{i}}=([(\alpha_i,\alpha_i)]_q)^{-n}
(\tilde{\rm ad}\; e_{\delta}^{(i)})^{n}
e_{\delta -\alpha_{i}} ,\qquad
e_{-(n+1)\delta +\alpha_{i}}=(-1)^n([(\alpha_i,\alpha_i)]_q)^{-n}
(\tilde{\rm ad}\; e_{-\delta}^{(i)})^{n}
e_{-\delta +\alpha_{i}},
$$
$$
{e'}_{(n+1)\delta}^{(i)}=
[e_{n\delta + \alpha_{i}}, e_{\delta - \alpha_{i}}]_{q}  ,\qquad
{e'}_{-(n+1)\delta}^{(i)}=
[ e_{-\delta + \alpha_{i}},e_{-n\delta - \alpha_{i}}]_{q^{-1}}
$$
(for $n > 0$), where $(\tilde{\rm ad}\,x)y=[x,y]$ is a usual
commutator.

The imaginary root vectors
$e^{(i)}_{\pm n \delta}$, which coincide with $a_{i,\pm n}$ up to
central factors,
 are related to ${e'}_{n\delta}^{(i)}$ via Schur polynomials, namely
$$E_{\pm i}(z)=  \ln (1 + E'_{\pm i}(z)),$$
where
$$E_{\pm i}(z)= \pm(q-q^{-1})\sum_{m\geq 1}^{} e_{\pm m\delta}^{(i)}z^m,
\qquad
E'_{\pm i}(z)= \pm(q-q^{-1})\sum_{m\geq 1}^{} {e'}_{\pm m\delta}^{(i)}z^m.$$

The rest of the real root vectors we construct in accordance with
the induction procedure \r{induction} using the root vectors $e_{n\delta +
 \alpha_i}$,
$e_{(n+1)\delta - \alpha_i}$, $e^{(i)}_{(n+1)\delta}$,
($i = 1,2, \ldots, r$; $n \in {\ZZ_{+}}$). In this procedure the dependence
 on the choice of the last vector, that is on the choice of an index $i$ is
 only
 in the normalization of the root vectors.

For any real root $\gamma$ denote by $R_\gamma$ the formal series
$$R_\gamma= \exp_{q^{(\gamma,\gamma)}}\left(C_{\gamma}^{-1}e_\gamma\ot
e_{-\gamma}\right)$$
where $C_{\gamma}$ is the normalization constant given from the relation
$[e_{\gamma},e_{-\gamma}]=C_{\gamma}\left(k_\gamma-k_\gamma^{-1}\right)$ and
$$\exp_p(x)=1+\sum_{n>0}\frac{x^n}{(n)_p!}\ ,
\qquad (n)_p=\frac{p^n-1}{p-1}\ .$$
Put also
\bn \K^{} =
q^{-t}q^{\frac{-c\ot d -d\ot c}{2}}
\prod_{n>0}\exp\left(-n(q-q^{-1})\sum_{i,j=1}^r d_{i,j}^{(n)}
a_{i,n}\ot a_{j,-n}\right)
q^{\frac{-c\ot d -d\ot c}{2}},
\label{Cartan}
\ed
where
 $t=\sum h_i\ot h^i$ is an invariant element in tensor square of
 Cartan subalgebra $\hh\ot\hh$ of $\g$; $d_{i,j}^{(n)}$ is an inverse matrix
to symmetrized Cartan matrix of $\g$
 $$b_{i,j}^{(n)}=[n(\a_i,\a_j)]_q\ , \qquad i,j=1,...,r\ .$$
The universal $R$-matrix $\kr^{CW}$ for quantum affine algebra $\Uqg$ with
Drinfeld comultiplication $\Delta^{(1)}$ was identified in \cite{KT1} with
formal ordered infinite product
 \bn
\kr^{CW}=\K\kr^{CW}_-\kr^{CW}_+\ ,
\label{RCW}
\ed
where
$$\kr^{CW}_-=\prod\limits^{\to}_{\gamma\in \underline{\Delta}^{re},\gamma>
\delta}
R_\gamma, \qquad
\kr^{CW}_+=\left(\prod\limits^{\to}_{\gamma\in
\underline{\Delta}^{re},\gamma>
\delta}
R_\gamma\right)^{21}.$$
 For instance, in the
case of  $\Uqdva$ the general expression \r{RCW} looks
as
\bn
\kr^{CW}=q^{-\frac{h\ot h}{2}}q^{-c\ot d-d\ot c}
\exp\sk{(q^{-1}-q)\sum_{k>0}\frac{n}{[2n]_q}\ q^{-\frac{cn}{2}} a_n\ot
q^{\frac{cn}{2}} a_{-n}}\cdot
\stackrel{\longrightarrow}{\prod_{n\in\z}}
\exp_{q^2}\Big((q^{-1}-q^{})\xm_{-n}\otimes \xp_{n}\Big)
\label{RCW1}
\ed
where we drop for simplicity the index of simple root
and put $e_n\equiv e_{\a,n}$, $f_n\equiv e_{-\a,n}$.
\begin{prop}
\label{action}
The tensor $\kr^{CW}$ from \r{RCW} defines correctly determined operator in
tensor product $W\ot V$ of $\Uqg$-modules if either $W$ is a lowest weight
representation or $V$ is a highest weight representation.
\end{prop}
Let us consider a case when $V$ is highest weight representation.
Take two weight vectors  $v\in V$ and $\xi\in V^*$. Consider the matrix
coefficient $\<\xi, \R^{CW} v\>$. We claim that only finitely many terms of
formal series \r{RCW} contribute into this matrix coefficient. Indeed, let
$\rho\in\hhh^*$, where $\hhh$ is Cartan subalgebra of $\ggg$, satisfy the
conditions $(\a_i,\rho)=1$, $i=0,1,..,r$. Then, by definition of highest
weight representation, there exists some big positive $N$ such that
 $xv=0$ for any $x\in \Uqg$, such that $(\la(x),\rho)>N$, where $\lambda(x)$
is a weight of $x$. All the entries from $\kr^{CW}_+$ have positive weight
with respect to right tensor factor and there are only finitely many terms
that satisfy to condition $(\la(x),\rho)<N$. Let us fix one such term $x$.
 We are restricted now to the summands of the type $yx$,
where $y\in \K\kr_-^{CW}$, so for the $y$ we have opposite restriction
$(\la(y),\rho)>-M$ for some positive $M$ because $\xi\in V^*$. But the
weights of the terms from $\K\kr_-^{CW}$ have negative weights with respect
to the right tensor factor, so by the same arguments we have only finitely
many choices for $y$.

This shows that any matrix coefficient $\<\eta\ot\xi, R w\ot v\>$,
 $v\in V, \xi\in V^*, w\in W, \eta\in W^*$ is well
define finite sum if $V$ is of highest weight. Analogous arguments are valid
when
$W$ is of lowest weight.

Proposition \ref{action} allows to check the identities for
$\kr^{CW}$ by their
application  onto highest weight or lowest
weight representation. For instance, the equalities
$$
\Delta^{op}(a)=\kr^{CW}\Delta(a)\left(\kr^{CW}\right)^{-1},
$$
$$
(\Delta\ot {\rm id})\kr^{CW}=\left(\kr^{CW}\right)_{13}\left(\kr^{CW}
\right)_{23},\qquad
( {\rm id}\ot \Delta)\kr^{CW}=\left(\kr^{CW}\right)_{13}\left(\kr^{CW}
\right)_{12}.
$$
are correct at least when we apply them to tensor products of highest (or
lowest) weight representations.

 \subsection{Analytical properties of matrix coefficients \cite{E}}

  Our considerations are strongly based on the analytical properties
 of the matrix coefficients of the products of the generating functions
 for $\Uqg$ \cite{E}. In the most general form they can be formulated as
follows. Let $a_k(z)$ stands for any of the generating functions of the type
$e_{\pm\a_i}(z)$,  $\psi^\pm_{\a_i}(z)$, $V$ be a highest degree
representation in a sense of the previous section, $v\in V$ and $\xi\in V^*$
 be two homogeneous vectors.
 Consider  the matrix coefficient
$$\<\,\xi, a_{k_1}(z_1)\cdots a_{k_m}(z_m) v\,\>$$
as a formal power series over $z_1,..., z_m$.

Then, first, this formal power series belongs to a space
\bn
 {\CC}[z_1,z_1^{-1},...,z_m,z_m^{-1}][[\frac{z_2}{z_1},
\frac{z_3}{z_2},...,\frac{z_m}{z_{m-1}}]],
 \label{zona}
 \ed
that is, can be presented as Taylor series over the variables
  $z_2/z_1,..., z_{m-1}/z_m$ with coefficients being polynomials over
 $z_1, z_1^{-1},...,z_m, z_m^{-1}$.

 Second, this formal power series converges in the region
 $|z_1|\gg|z_2|\gg...\gg|z_m|$ to a rational function which poles in
 $\left(\CC^*\right)^m$ are dictated by the quadratic
 relations \r{1}--\r{10}.

 The proof of the first statement is based on the observation that under the
 conditions on the vectors $v$ and $\xi$ for some $N$ big enough the
coefficient
at $z_1^{-n_1}z_2^{-n_2}\cdots z_m^{-n_m}$ vanishes if $n_m>N$,
or $n_{m-1}+n_{m}>N$, or
  $n_{m-2}+n_{m-1}+n_{m}>N$  and so on due to the
definition of highest degree representation.

 Once \r{zona} is established, we can repeat the arguments of
 \cite{E} in order to show that the formal power series converges to a
 meromorphic function.

 Consider, for instance, the matrix coefficient
  $\<\,\xi, e_{\a_{i_1}}(z_1)\cdots e_{\a_{i_m}}(z_m )v\,\>$.
 We know from commutation relations that
 \bn
\prod_{k<l}(z_k-q^{(\a_{i_k},\a_{i_l})}z_l)
 e_{\a_{i_1}}(z_1)\cdots e_{\a_{i_m}}(z_m )=
 \prod_{k<l}(q^{(\a_{i_k},\a_{i_l})}z_k-z_l)
 e_{\a_{i_m}}(z_m)\cdots e_{\a_{i_1}}(z_1 ),
\label{mer1}
\ed
 so the matrix coefficient
  $$\<\,\xi,
 \prod_{k<l}(z_k-q^{(\a_{i_k},\a_{i_l})}z_l)
e_{\a_{i_1}}(z_1)\cdots e_{\a_{i_m}}(z_m )v\,\>$$
  belongs  to
 $\CC[z_1,z_1^{-1},...,z_m,z_m^{-1}][[\frac{z_2}{z_1},
\frac{z_3}{z_2},...,\frac{z_m}{z_{m-1}}]]$; but due to the form of the
r.h.s. of
\r{mer1}, it belongs also to
 $\CC[z_1,z_1^{-1},...,z_m,z_m^{-1}][[\frac{z_1}{z_2},
\frac{z_2}{z_3},...,\frac{z_{m-1}}{z_{m}}]]$,
 so lies in  their intersection,
 $\CC[z_1,z_1^{-1},...,z_m,z_m^{-1}]$
 It means that  the series
  $\<\xi, e_{\a_{i_1}}(z_1)\cdots e_{\a_{i_m}}(z_m )v\>$.
 converges in the region $|z_k|>|q^{(\a_{i_k},\a_{i_l})}z_l|$
  to a meromorphic in $\left(\CC^*\right)^m$ function with simple
  poles at $z_k=q^{(\a_{i_k},\a_{i_l})}z_l$.

 These properties of matrix coefficients allow to treat in the highest
 weight  representation the generating functions
 $a_{k_1}(z_1)\cdots a_{k_m}(z_m)$
as operator valued functions, analytical in
 the region  $|z_1|\gg|z_2|\gg...\gg|z_m|$, where it coincides with the
 products
 $\left(a_{k_1}(z_1)\cdots a_{k_l}(z_l)\right)\cdot
 \left(a_{k_{l+1}}(z_{l+1})\cdots a_{k_m}(z_m)\right)$ for any
 $l: 1\leq l\leq m$. This is because the multiplication
 of Taylor series is well defined. Moreover, the relations \r{1}--\r{10}
 describe the
 analytical continuations of these functions to other regions. For instance,
  the analytical continuation of the function
 $e_{\al_i}(z)e_{\al_j}(w)$  from the region
 $|z|>|q^{({\al_i},{\al_j})}w|$ to the region
 $|z|<|q^{({\al_i},{\al_j})}w|$ is defined, due to \r{1}, by formal power
 series in the completed algebra $\UqgD$
 $$\frac{q^{-({\al_i},{\al_j})}-\frac{z}{w}}
{1-q^{-({\al_i},{\al_j})}\frac{z}{w}}
 e_{\al_j}(w)e_{\al_i}(z).$$

\subsection{The configuration spaces}

We are interested now in the operator valued analytical functions, given as
products (in the above sense) of the currents
$$t_{\a_i}(z)=(q^{-1}_{}-q_{})e_{-\a_i}(z)\otimes
e_{\a_i}(z).$$
We have, due to the relations \r{1}:

The function $t_{\a_{i_1}}(z_1)\cdots t_{\a_{i_n}}(z_n)$ is analytical
 in the region
$|z_k|>{\rm max} (|q^{\pm(\a_{i_k},\a_{i_l})}z_l|$, $k<l$, admits symmetric
 meromorphic
analytical continuation to $\left(\CC^*\right)^n$ with simple poles at
shifted diagonals $z_k= q^{\pm(\a_{i_k},\a_{i_l})}z_l$. Here symmetric group
exchanges simultaneously the variables and their root labels. For instance,
for $n=2$ it means that the function $t_\a(z)t_\b(w)$ has two simple
 poles at
$z=q^{(\a,\b)}w$ if $(\a,\b)\not=0$ and is commutative
$$t_\a(z)t_\b(w)=t_\b(w)t_\a(z)$$
in a sense of analytical continuation.

Let now $\Pi$ be a set of simple positive roots of $\g$ and $I$ be a finite
 set
$k_1,k_2,...,k_n$ of integers equipped with labels of simple roots, that is,
$I$ is a finite subset $\II=\{k_1,k_2,...,k_n\}\subset {\NN}$ of a set of
positive integers together with a map $\iota:\II \to \Pi$.
Let $X_I$ be the following stratified space. As a total space, $X_I$ is
isomorphic to ${\CC}^n$ with coordinates $z_k, k\in \II$. The closures
of the strata are
given by the intersections of hyperplanes
 $H_{k,l}=\{z_k=q^{(\io(k),\io(l))}z_l\}$ for
any $k,l\in\II$ such that $(\io(k),\io(l))\not=0$ and
 $H_i=\{z_{i}=0, i\in\II\}$.
 By $U_I$ we denote an open stratum: the complement to the union of
hyperplanes.

Among all the strata there are the distinguished ones, which we call Serre
strata. Let $\a$ and $\b$ be two simple adjacent roots, $a_{\a,\b}$ be
 corresponding
element of Cartan matrix, $m_{\a,\b}=1-a_{\a,\b}$. Let $z_{l_1},\ldots,
 z_{l_m}$
 be the coordinates of $X_I$, labeled by $\a$, that is $\io(l_i)=\a,
i=1,\ldots,m$ and $z_{l_{0}}$ be a coordinate labeled by $\b$:
$\io(l_{0})=\b$. Then the stratum
\bn
z_{l_{0}}=q^{-\frac{m(\a,\a)}{2}}z_{l_1}=q^{-\frac{(m-2)(\a,\a)}{2}}z_{l_2}
=\ldots=
q^{\frac{(m-2)(\a,\a)}{2}}z_{l_{m-1}}=q^{\frac{m(\a,\a)}{2}}z_{l_m}
\label{Serst1}
\ed
is called Serre stratum. Another type of Serre strata appear due to the
 vanishing of the squares of the same
 fields in the same point. They have a form
\bn
z_{m_{1}}=z_{m_{2}}=q^{\pm (\io(m_1),\io(m_3))}z_{m_3}
\label{Serst2}
\ed
 for any $m_1, m_2, m_3$ such that $\io(m_1)=\io(m_2)$.

We are interested in the integrals of the $n$-forms witch have simple
poles
 at the hyperplanes $H_{k,l}$ and meromorphic singularities at
hyperplanes $H_k$ over certain $n$-cycles
 in $U_I$:
\bn
\omega=\frac{P(z,z^{-1})}
{(2\pi i)^n\prod_{l\not=m}(z_{k_l}-q^{(\io(k_l),\io(k_m))}z_{k_m})}
\frac{dz_{k_1}}{z_{k_1}}\wedge \frac{dz_{k_2}}{z_{k_2}}\wedge ...\wedge
\frac{dz_{k_n}}{z_{k_n}}
\label{forma}
\ed
where $P(z,z^{-1})$ is a polynomial over $z_{k_i}$ and $z_{k_i}^{-1}$,
$n=|\II|$.

The description of the homologies of the complement to the arrangement of
hyperplanes is well known, see, e.g. \cite{Varchenko}.
It looks as follows. Suppose we have a collection of hyperplahes $H_i$
 in a space $X={\CC}^N$. Hyperplanes define on $X$ the structure of
stratified space; the (closed) strata are all possible intersections of
 hyperplanes. A sequence
$$L_n\subset L_{n-1}\subset ... \subset L_0$$
of strata is called a (full) flag of length $n$, if ${\rm codim}(L_k)=k$

Let $U$ be a complement to the union of hyperplanes.
 The  space $H_n(U,{\CC})$ is isomorphic
 to a factor of vectorspace, generated by all flags of length $n$ modulo
 the Orlik-Solomon relations, attached to any incomplete flag
 $L_n\subset ...\subset L_{k+1}\subset L_{k-1}...\subset L_0,$
${\rm codim}( L_i)=i$. The relation is
$$\sum_{{L_k:L_{k+1}\subset L_k\subset L_{k-1}\atop {\rm codim} L_k=k}}
L_n\subset ...\subset L_{k+1}\subset L_k\subset L_{k-1}...\subset L_0=0.$$

The cycle, corresponding to a flag
$L_n\subset L_{n-1}\subset ... \subset L_0$
is $n$-dimensional torus, which is described as follows:
first we take a small circle with a center in a generic point of $L_n$,
 surrounding $L_n$ inside $L_{n-1}$, then for
 any point of this circle we draw a circle surrounding $L_{n-1}$ inside
$L_{n-2}$ such that it does not intersect other strata inside $L_{n-1}$ and
 so on. Note that in this procedure the next circle is  much smaller
 then the previous one. The resulting torus lies in $U$ and its orientation
is given by the order in the procedure above.

For any cycle
$L=\{L_n\subset L_{n-1}\subset ... \subset L_0\}$
and for any $m$-form $\omega$, regular in $U$ we denote by
 ${\rm Res}_L \omega$ the integral of the form $\omega$ over the cycles,
defined
by $L$ (we take into account their dependence over generic point of $L_n$).
It is an $(m-n)$-form on the open part of $L_n$.

Return now to stratified space $X_I$.
Denote by $\Omega_I$ the space of such forms \r{forma}, for which any
repeated residue to any Serre stratum \r{Serst1}, \r{Serst2} is zero.
Namely, we say that $\omega\in\Omega_I$ if ${\rm Res}_L \omega=0$ for any
$L=\{L_n\subset L_{n-1}\subset ... \subset L_0\}$, such that $L_n$ is a
Serre stratum \r{Serst1}, \r{Serst2}. We call such forms admissible.

The main goal of this subsection is the definition
of factorizable  system of antisymmetric cycles
for any  quantum affine algebra $\Uqg$ .
 This definition means an assignment to any labeled
 set $I$ a
  symmetric $n$-cycle $D_I\in H_n(U_I)$, where $n=|I|$, with certain
factorization properties. Let us first explain symmetricity condition.

Any permutation $\sigma\in S_n$ where $n=|I|$ defines a new ordered
set $\sigma(I)$ with induced labeling. Moreover, $\sigma$ defines a
diffeomorphism of configuration spaces: $\sigma:X_I\to X_{\sigma(I)}$.
 We demand that:
\bn
\oint\limits_{D_{\sigma(I)}} \omega=\oint\limits_{D_I}\sigma^*(\omega)
\label{symm}
\ed
for any $\sigma\in S_{|I|}$ and $\omega\in \Omega_{\sigma(I)}$.
Equivalently,
$$\oint\limits_{D_{\sigma(I)}}
\omega=(-1)^{l(\sigma)}\oint\limits_{\sigma(D_I)}\omega,$$
which means that the cycle $D_I$ is defined on the factor of $U_I$ over the
action of the product $S_{i_1}\times S_{i_2}\times ...\times S_{i_r}$
of symmetric groups, where $i_k=\# k\in\II, \iota(k)=\a_k$ and can
be labeled by
an element of the lattice of positive roots of $\g$.

Suppose now that finite labeled set $I$, $|I|=n+m$ is presented as a disjoint
union of its labeled subsets $I=I^1\coprod I^2$, $|I^1|=n, |I^2|=m$. Then for
 any $D_{I^1}\in H_n(U_{I^1}, {\CC})$ and for any $D_{I_2}\in
H_m(U_{I^2}, {\CC})$ we can
 attach the cycle $D_{I^1}\times D_{I^2}\in H_{n+m}(U_I, {\bf C})$ as
follows.
There is
a natural map $\phi:X_{I^1}\times X_{I^2}\to X_I$. The cycle $D_{I^2}$ is
$m$-dimensional compact manifold, and we can apply to it arbitrary
dilatation
$z_l\to \varepsilon z_l, l\in \II^2, \varepsilon>0$ without changing
 homology
class. We choose such small $\varepsilon$ that all the points of the direct
product of dilated $D_{I^1}$ and of $D_{I^2}$ belong to $U_{I}$ under the
image of $\phi$. It is possible by the following reasons: on every
$D_{I^j}$ absolute values $|z_k|$ of the coordinates are restricted from
both
sides, which follows from the compactness of the cycles and from the
observations that coordinate hyperplanes $z_k=0$ should not intersect the
cycles by their definition. Then it is clear that for small enough
 $\varepsilon$ a direct
product of dilated $D_{I^1}$ and of $D_{I^2}$
does not intersect hyperplanes $z_k=q^{\pm(\io(k),\io(l))}z_l$ ,
where $k\in \II^1$ and $l\in \II^2$.
But this is the only thing we are checking.
This construction  defines invariantly the cycle $D_{I^1}\times D_{I^2}\in
H_{n+m}(U_I, {\CC})$ equipped with natural orientation.

 Let us choose any form  $\omega\in \Omega_I$
satisfying additional property: it has no singularities at hyperplanes
$\{z_k=q^{(\io(k),\io(l))}z_l\}$ for all $k\in I^1$ and $l\in I^2$.
Denote the space of such $\omega$ by $\Omega_{I^1,I^2}$.
We demand that the integral of any $\omega\in\Omega_{I^1,I^2}$
over $D_{I}$ coincides with the
integral over $D_{I^1}\times D_{I^2}$:
\bn
\oint\limits_{D_I}\omega=\oint\limits_{D_{I^1}\times D_{I^2}}\omega
\label{razbienie}
\ed
The factorization conditions mean that, first, the relation \r{razbienie}
holds for any decompositions $I=I^1\coprod I^2$
and  for any $\omega\in \Omega_{I^1,I^2}$ and,
second, the initial conditions
\bn
D_{\{k\}}= \{|z_k|=1\}
\label{initial}
\ed
take place for any one point set $I$, $\II=\{k\}$.

 Let $q^n\not=1$ for any natural $n$. We suggest the following
\bigskip

\noindent{\bf Conjecture}\ \ \
{\it
  For any $\Uqg$, where $\g$ is simple Lie algebra, there exists a factorizable
system $\DD_I$
of antisymmetric cycles with initial conditions \r{initial}.}
\bigskip

\noindent
Moreover, we suppose that the cycles are unique as the functionals over
admissible forms. Also, we suppose that they are uniquely defined by
factorization conditions; the symmetricity should follow from the
factorization properties.

 We can describe precisely such a system for the algebra $\UqdvaD$.
In this case  we have no Serre restrictions \r{Serst1} and, since there
is only
one simple root for Lie algebra $\g$, the cycles $\DD_I$ are
 parameterized by
positive integers $n$, so we denote them as $D_n$. They look as follows.
  Fix some order of the variables $z_{k_1}, z_{k_2},...,z_{k_n}$.
Then in the flag description
the cycle $D_n$ consists of all $n$- flags
 $L=\{L_n\subset L_{n-1}\subset ... \subset
L_0\}$ in $X_n={\CC}^n$, all with multiplicity one, such that $L_1$ is one
of the
hyperplanes $$z_{k_1}=0,\quad z_{k_1}=q^{-2}z_{k_j}\quad,
 j\not=1,$$ $L_2$ is codimension two stratum, which can be obtained as an
 intersection of one of hyperplanes $L_1$, listed above, with hyperplanes
 $$z_{k_2}=0,\quad z_{k_2}=q^{-2}z_{k_j},\quad j\not=2$$
 and so on. Clear, that in the last stage $L_n$ coincides with an origin
 in ${\CC}^n$.

\begin{prop}
\label{sldvatheorem}
The cycles $D_n$ form factorizable system of antisymmetric cycles for
$\Uqdva$, if $q$ is not root of $1$.
\end{prop}

There are two short proofs of the proposition.
Note first that it is sufficient to prove the symmetricity condition only.
Then the factorization property follows from the definition of the cycles
$\DD_I$.

For the first proof, it sufficient to check that for $|q|>1$, the cycles
$D_n$
 for $I=\{1,...,n\}$ are homotopic to the product of unit circles, which are
 clearly antisymmetric  and the factorization property is also clear.
For another proof, one can note that the integral of $\omega$ over $D_n$
coincides with Grothendick residue \cite{GH} of $\omega$ with respect to the
system of divisors
\bn
f_i(z)=z_i\prod_{j\not=i}(z_i-q^{-2}z_j)=0, \qquad i=1,2,...,n.
\label{divisor}
\ed
 The Grothendique residue is correctly defined for a system of divisors
$f_i(z)=0$ with common intersection being a point $\{z_i=0\}$ and can be
 written in this case as an integral
$\oint\limits_\Gamma \omega$, where
 $\Gamma=\{z:|f_i(z)|=\varepsilon\}$. The system \r{divisor} satisfy these
conditions. By the definition of the residue,
it is  antisymmetric. The factorization property follows then from the form
of $f_i(z)$.

Note that the vanishing conditions on 'diagonal' Serre strata \r{Serst2}
were not used in the proof.
\subsection{The properties of the universal $R$-matrix and factorization of
 the cycles}
The universal $R$-matrix $R$ for a quasitriangular Hopf algebra ${\cal A}$ is
characterized
 by the properties
\bn
\Delta^{op}(a)=R\Delta(a)R^{-1}
\label{intertwining}
\ed
 for any $a\in {\cal A}$ and
\bn
(\Delta\ot {\rm id})R=R_{13}R_{23},\qquad
( {\rm id}\ot \Delta)R=R_{13}R_{12}.
\label{quasi}
\ed
Let now ${\cal A} =\UqgD$ where ${\frak g}$ has rank $r$
and $D_I$ be a factorizable system of
antisymmetric cycles. Put
\bn {\cal R}=\K\R \label{RR} \ed where
$\K$ is given in \r{Cartan}
and
$$
\R=\Pqexp\oint \dz{}\sum\limits_{i=1,...,r}t_{\a_i}(z)=$$
\bn
1+\sum_{n>0}\frac{1}{n!}\sum_{I:\II=\{1,2,...,n\}}\;\oint\limits_{D_I}
 t_{\io(1)}(z_1)t_{\io(2)}(z_2)\cdots t_{\io(n)}(z_n)\dzz{1}\wedge\cdots\wedge
\dzz{n}
\label{Pexp}
\ed
Here and below $\dzz{k}$ means $\frac{dz_k}{2\pi iz_k}$.
\begin{th}
\label{theorem1}
For any factorizable system $\{D_I\}$ of antisymmetric cycles the tensor
 \r{RR} satisfy the properties \r{intertwining} and \r{quasi}
of the universal
$R$-matrix for topological Hopf algebra $\UqgD$
with comultiplication $\D^{(1)}$.
\end{th}
 Note first that the form
 $t_{\io(1)}(z_1)t_{\io(2)}(z_2)\cdots t_{\io(n)}(z_n)\dzz{1}\wedge
\cdots\wedge
\dzz{n}$
  satisfy vanishing conditions on the Serre strata. It can be deduced
 from the Serre relations in a form, presented
in \cite{E}. We omit corresponding calculations. By
 antisymmetricity of the cycles we can rewrite then the series \r{Pexp}
 in the form
$$\R=\sum_{n_1,...,n_r\geq  0}\frac{1}{n_1!\cdots n_r!}
\oint_{D_{n_1,...,n_r}}
 t_{\io(1)}(z_1)t_{\io(2)}(z_2)\cdots t_{\io(n)}(z_n)\dzz{1}\wedge
\cdots\wedge
\dzz{n}$$
where $D_{\{n_i\}}$ is the cycle, corresponding to labeled set
$I$, $\II=\{1,...,n_1+\cdots +n_r\}$ and $\io(k)=\a_i$ if
$n_1+\cdots +n_{i-1}<k\leq n_1+\cdots +n_i$ for any order $\a_1, ... ,\a_r$
 of simple roots of ${\frak g}$.

Let us prove an equality \r{intertwining} for $x=e_{\a_1}(z)$.
The basic property of the tensor $\K$, which we use
here, is
\bn
\D^{(1),op}(x)=\K\D^{(2)}(x)\K^{-1}
\label{K}
\ed
for all $x\in \UqgD$.
We have ($\Delta$  means $\Delta^{(1)}$):
$$\kr \Delta e_{\a_1}(z) \kr^{-1} =
 \K \R(e_{\a_1}(z)\otimes 1)( \R)^{-1}\K^{-1}  +
\K \R(\psi^-_{\a_1}(zq^{\frac{c_1}{2}})\otimes e_{\a_1}(zq^{c_1}) )
( \R)^{-1}\K^{-1}. $$
Denote the summands of r.h.s. as
$X_1=\K \R(e_{\a_1}(z)\otimes 1)( \R)^{-1}\K^{-1}$
and $X_2= \K \R(\psi^-_{\a_1}(zq^{\frac{c_1}{2}})\otimes e_{\a_1}(zq^{c_1}) )
(\R)^{-1}\K^{-1} $.
We see from \r{4} and \r{10} that  $(n=n_1+...+n_r)$
$$X_1-
\K (e_{\a_1}(z)\otimes 1)\K^{-1}=
\K\sum_{{n_2,...,n_r\geq 0\atop n_1\geq k>0} }
\frac{1}{n_1!\cdots
n_r!} \oint\limits_{D_{\{n_i\}}}\dzz{1}\wedge...\wedge\dzz{n}
t_{\a_1}(z_1)\cdots t_{\a_1}(z_{k-1})$$
$$\left(\left( \delta(\frac z {z_k}q^{-c})\psi^+_{\a_1}(z_kq^{\frac{1}{2}c})-
          \delta(\frac z {z_k}q^{c})\psi^-_{\a_1}(zq^{\frac{1}{2}c})
\right)\otimes e_{\a_1}(z_i)\right)
t_{\a_1}(z_{k+1})\cdots
t_{\a_1}(z_{n_1})
\prod\limits_{\io(j)>1}t_{\io(j)}(z_j)
  \R^{-1}\K^{-1}= $$
$$
\K\sum_{{n_2,...,n_r\geq 0\atop n_1>0} }
\frac{n_1}{n_1!\cdots
n_r!} \oint\limits_{D_{\{n_i\}}}
 \delta(\frac{z}{z_1q^{c}})\psi^+_{\a_1}(z_1q^{\frac{1}{2}c})
\otimes e_{\a_1}(z_1)
t_{\a_1}(z_2)\cdots t_{\a_1}(z_{n_1})
\prod\limits_{\io(j)>1}t_{\io(j)}(z_j)\dzz{1}\wedge...\wedge\dzz{n}
  \R^{-1}\K^{-1} $$
$$ -
\K\sum_{{n_2,...,n_r\geq 0\atop n_1>0} }
\frac{n_1}{n_1!\cdots
n_r!} \oint\limits_{D_{\{n_i\}}}
t_{\a_1}(z_2)\cdots t_{\a_1}(z_{n_1})
\prod\limits_{\io(j)>1}t_{\io(j)}(z_j)
 \delta(\frac {zq^c}{z_1} )\psi^-_{\a_1}(zq^{\frac{1}{2}c})
\otimes e_{\a_1}(z_1)
\dzz{1}\wedge...\wedge\dzz{n}
  \R^{-1}\K^{-1} $$
 due to antisymmetricity of the cycle $D_{\{n_i\}}$.
    {}From the commutation relations \r{1}, \r{3}, \r{4} it follows that
 the integrands  belong
 to a class $\Omega_{I_{\{n_i\}}}$ and, moreover, the first integrand has no
singularities at $z_j=q^{(\a_1, \io(j))}z_1$ for all $j\not=1$ and the second
integrand has no singularities at
$z_1=q^{(\a_1, \io(j))}z_j$ for all $j\not=1$.
{}From factorization condition and from the definition of the multiplication
 of the currents together with initial condition \r{initial} we see that
$$X_1-
\K(e_{\a_1}(z)\otimes 1)\K^{-1}= $$
$$
\K
\oint\dzz{1}
 \delta(\frac{z}{z_1q^{c}})
\psi^+_{\a_1}(z_iq^{\frac{1}{2}c})
\otimes e_{\a_1}(z_1)\cdot
\sum_{{n_2,...,n_r\geq 0\atop n_1>0}}
\frac{1}{(n_1-1)!\cdots n_r!}
 \oint\limits_{D_{n_1-1,n_2,...,n_r}}
\prod\limits_{i>1}t_{\io(i)}(z_i)\dzz{2}\wedge...\wedge\dzz{n}
  \R^{-1}\K^{-1} $$
$$-
\K
\sum_{{n_2,...,n_r\geq 0\atop n_1>0}}
\frac{1}{(n_1-1)!\cdots n_r!}
 \oint\limits_{D_{n_1-1,n_2,...,n_r}}
\prod\limits_{i>1}t_{\io(i)}(z_i)\dzz{2}\wedge...\wedge\dzz{n}
  \cdot
\oint\dzz{1}
 \delta(\frac{zq^c}{z_1})\psi^-_{\a_1}(zq^{\frac{1}{2}c})
\otimes e_{\a_1}(z_1)\R^{-1}\K^{-1}$$
$$
=\K\psi^+_{\a_1}(zq^{-\frac{c_1}{2}})
\otimes e_{\a_1}(zq^{-c_1})\K^{-1}-
\K\R
\psi^-_{\a_1}(zq^{\frac{c_1}{2}})
\otimes e_{\a_1}(zq^{c_1})\R^{-1}\K^{-1}$$
after a change of the index of summation over $n_1$. The last summand cancels
 $X_2$ so
$$\kr \Delta e_{\a_1}(z) \kr =
\K(e_{\a_1}(z)\otimes 1)\K^{-1}+
\K\psi^+_{\a_1}(zq^{-\frac{c_1}{2}})
\otimes e_{\a_1}(zq^{-c_1})\K^{-1}=
\K\Delta^{(2)}(e_{\a_1}(z))\K^{-1}=
\Delta^{(1),op}(e_{\a_1}(z))$$
 which finishes the proof.

Let us prove one of the relations \r{quasi}.
The left hand side equals to
$$(\Delta\ot 1)\K\R= \K_{13}\K_{23}(\Delta\ot 1)\R=
\K_{13}\K_{23}
\Pqexp\oint \dz{}\sum\limits_{i=1,...,r} \left( t^1_{\a_i}(z)+
t^2_{\a_i}(z)\right),$$
where
\bn
 t^1_{\a_i}(z)=(q^{-1}_{}-q_{})e_{-\a_i}(zq^{c_2})\ot\psi^+_{\a_i}
(zq^{c_2/2})\ot e_{\a_i}(z)
\label{t1}
\ed
and
\bn
 t^2_{\a_i}(z)=(q^{-1}_{}-q_{})\cdot 1\ot e_{-\a_i}(z)\ot e_{\a_i}(z)
\label{t2}
\ed
On the other side,
$$\kr_{13}\kr_{23}=\K_{13}\R_{13}\K_{23}\R_{23}=\K_{13}\K_{23}
\left(\K_{23}^{-1}\R_{13}\K_{23}\right)\R_{23}=$$
$$
\K_{13}\K_{23}\Pqexp\oint \dz{}\sum\limits_{i=1,...,r} t^1_{\a_i}(z)\cdot
\Pqexp\oint \dz{}\sum\limits_{i=1,...,r} t^2_{\a_i}(z)$$
 due to the properties of $\K$. So we have to prove an equality
\bn
\Pqexp\oint \dz{}\sum\limits_{i=1,...,r} \left( t^1_{\a_i}(z)+
t^2_{\a_i}(z)\right)=\Pqexp\oint \dz{}\sum\limits_{i=1,...,r} t^1_{\a_i}(z)
\cdot
\Pqexp\oint \dz{}\sum\limits_{i=1,...,r} t^2_{\a_i}(z)
\label{addition}
\ed
where $t^1_{\a_i}(z)$ and $t^2_{\a_i}(z)$ are given by \r{t1} and \r{t2}.
An equality \r{addition} is equivalent to following equalities of
integrals. Let $I$, $\II=\{1,...,n\}$ be a labeled set and $j(1),...,j(n)$ be
arbitrary sequence of numbers $1$ and $2$, $j(k)=1,2$. Then we should have:
$$
\oint\limits_{D_I}
 t_{\io(1)}^{j(1)}(z_1)t_{\io(2)}^{j(2)}(z_2)\cdots t_{\io(n)}^{j(n)}(z_n)
\dzz{1}\wedge\cdots\wedge
\dzz{n}=$$
\bn
\oint\limits_{D_I^1}
 t_{\io(k_1)}^{1}(z_{k_1})\cdots t_{\io(k_{m})}^{1}(z_{k_m})\dzz{k_1}\wedge
\cdots\wedge
\dzz{k_m} \cdot
\oint\limits_{D_I^2}
 t_{\io(k_{m+1})}^{2}(z_{k_{m+1}})\cdots t_{\io(k_{n})}^{2}(z_{k_{n}})
\dzz{k_{m+1}}\wedge\cdots\wedge
\dzz{k_n}
\label{addition2}
\ed
where $I^1=\{k_1,...,k_{m}; j(k_i)=1\}$ and $I^2=\{k_{m+1},...,k_{n};
j(k_i)=2\}$. Again, we see from commutation relations \r{7}, that the form
$\omega=t_{\io(1)}^{j(1)}(z_1)t_{\io(2)}^{j(2)}(z_2)\cdots t_{\io(n)}^{j(n)}
(z_n)\dzz{1}\wedge\cdots\wedge
\dzz{n}$ belongs to a class $\Omega_I$, and, moreover, has no singularities
at hyperplanes
$\{z_k=q^{(\io(k),\io(l))}z_l\}$ for all $k\in I^2$ and $l\in I^1$ and the
equality \r{addition2} holds due to factorization conditions on the cycles.

The two natural questions immediately appear for the $R$-matrix \r{RR}.
First, to which representations it can be applied, and second, how it is
connected to the $R$-matrix \r{RCW}, constructed by means of Cartan-Weyl
generators. Let us first try to apply \r{RR} to a tensor product $W\ot V$
 of the graded
representations with highest degree. Then there is a correct action of $\K$
and of each term in the series \r{Pexp} (for instance, because the integral
can be calculated by taking multiple residues which  belong to completed
algebra $\UqgD$). We can conclude now that this action is well defined in a
formal sense, that is, if we consider $\UqgD$ as an algebra over
$\CC[[q-1]]$ as well as the representations.

Let now $W$ and $V$ be highest weight representations. In particular, it
means that they are graded modules of highest degree. In this case we know
from Proposition \ref{action}
that the $R$-matrix \r{RCW} is well defined operator in $W\ot V$. In
particular, it means that $\kr^{CW}$ is also well defined operator in a
formal sense, that is, over the ring $\CC[[q-1]]$. By uniqueness argument we
conclude that in a formal sense the actions of $\kr^{CW}$ and of $\kr$
coincide. So they coincide as formal power series over $(q-1)$, and the
series \r{Pexp} converge and its action coincide with the action of
$\kr^{CW}$.
We summarize it in the following proposition.
\begin{prop}
The action of $R$-matrix $\kr$ \r{RR} on tensor product of highest weight
representations is well defined and coincides with the action of $\kr^{CW}$
 \r{RCW}.
\end{prop}

\section{Calculations for $\UqdvaD$}

In this section   we drop everywhere an
index of a simple root of $\sldva$ and denote
$e(z)\equiv e_\al(z)$ and $f(z)\equiv e_{-\al}(z)$.
 Following \cite{DK} and general receipt
of vertex
operator algebras, we study first the fields, appearing as the poles in the
products of generating functions. Surprisingly, we will see in the first
subsection, that the algebra, generated by the descenders of the tensor
field $f(z)\ot e(z)$,
 has simple closed form \r{T}. Next, we derive from the contour integral
description the recurrence relations on the terms of the universal
$R$-matrix for $\UqdvaD$, which can be written as simple differential
equation. Its solution has a structure of certain vertex operator due to the
simple structure of tensor fields. Finally, we present some application of
the derived expression for the $R$-matrix.

\subsection{Algebra of tensor fields}
\label{4.3}

 Let us first look to the relation \r{1}, which
now reads as $(z_1-q^2z_2)e(z_1)e(z_2)=(q^2z_1-z_2)e(z_2)e(z_1)$. It shows
 that the function
$e(z_1)e(z_2)$ has the only pole at the point $z_1=q^2z_2$. The residue
could be calculated in two ways. Since the product $e(z_1)e(z_2)$ is well
defined in completed algebra $\UqdvaD$, the contour of integration can be
replaced by two circles and, how it is explained in \cite{DK}, it gives the
residue as a difference of formal integrals:
$$ \res{z_1=q^2z_2}e(z_1)e(z_2)\ddz{1}= \oint
e(z_1)e(z_2)\dzz{1}-\oint\frac{q^{-2}-z_1/z_2}{1-q^{-2}z_1/z_2}e(z_2)e(z_1)
\dzz{1}
$$
which gives
\be\label{n8}
\res{z=q^2w}e(z)e(w)\ddz{}=
\sum_{n\in\z}w^{-2n}q^{-2n}\left\{(1-q^{-2})\xp_n^2+(q^2-q^{-2})
\sum_{k=1}^\infty q^{-2k}\xp_{n-k}\xp_{n+k}\right\}+\\ +\,
(q^2-q^{-2})\sum_{n\in\z}w^{-2n-1}q^{-2n-2}\sum_{k=0}^\infty
q^{-2k}\xp_{n-k}\xp_{n+k+1}
\ee
which is well defined series in completed algebra $\UqdvaD$.
We can use also the rule for the calculation of poles of the first order:
\be\label{rdef}
\res{z=q^2w}e(z)e(w)\ddz{}=
\lmt{z\to q^2w}\sk{1-q^2w/z}e(z)e(w)=(q^2-q^{-2})e(w)e(q^2w)
\ee
 Again, we can see from commutation relations that the product
 $e(z)e(w)e(q^2w)$ has unique pole, which equals, up to a constant, to
$e(w)e(q^2w)e(q^4w)$. Let us define the following generating functions from
 $\UqdvaD$
$e^{(n)}(z)$ and $f^{(n)}(z)$ by induction as
$$
e^{(n)}(z)=\res{z_1=q^{2(n-1)}z}e(z_1)e^{(n-1)}(z)\ddz{1},\qquad
f^{(n)}(z)=\res{z_1=q^{2(n-1)}z}f^{(n-1)}(z)f(z_1)\ddz{1},\qquad
$$
Put also
\bn
t^{(n)}(z)=-\res{z_1=q^{2(n-1)}z}t(z_1)t^{(n-1)}(z)\ddz{1},
\label{eft1}
\ed
where, as before, $t(z)=(q^{-1}-q)f(z)\ot e(z)$. and
$e^{(1)}(z)$ stands for $e(z)$,
$f^{(1)}(z)$  for $f(z)$,
$t^{(1)}(z)$  for $t(z)$.
Then, analogously to \r{rdef},
\bn
e^{(n)}(z)=(q^{-1}-q)^{n-1}[n-1]_q![n]_q!\tilde{e}^{(n)}(z),\qquad
f^{(n)}(z)=(q^{-1}-q)^{n-1}[n-1]_q![n]_q!\tilde{f}^{(n)}(z),
\label{eft}
\ed
\bn
\T{n}(z)=(q^{-1}-q)^{2n-1}[n-1]_q![n]_q!
\tilde{f}^{(n)}(z)\ot \tilde{e}^{(n)}(z)\ ,
\label{comp1}
\ed
 where
$$
\tilde{e}^{(n)}(z)=e(z)e(q^2z)\cdots e(q^{2(n-1)}z),\qquad
\tilde{f}^{(n)}(z)=f(q^{2(n-1)}z)\cdots f(q^{2}z)f(z)\ .
$$

Iterating the calculations \r{n8}, one can describe these fields in a
component form:
\bn\label{n18}
\tilde{e}^{(n)}(z)=\sum_{m\in\ZZ}\left(zq^{2n}\right)^m
\sum_{{\la_1\geq \cdots\geq \la_n,\atop \la_1+...+\la_n=m}}
\frac{  q^{-2(\la_1+...+k\la_k+...n\la_n)}}
{\prod_{j\in\ZZ}(\la'_j-\la'_{j+1})_{q^{2}}!}
e_{\la_n}e_{\la_{n-1}}\cdots e_{\la_1}
\ed
\bn\label{n19}
\tilde{f}^{(n)}(z)=
\sum_{m\in\ZZ}\left(zq^{-2}\right)^m
\sum_{{\la_1\geq \cdots\geq \la_n,\atop \la_1+...+\la_n=m}}
\frac{  q^{-2(\la_1+...+k\la_k+...n\la_n)}}
{\prod_{j\in\ZZ}(\la'_j-\la'_{j+1})_{q^{-2}}!}
 f_{\la_n}f_{\la_{n-1}}\cdots
f_{\la_1}
\ed
Here
 $\la'_j=\#
k$, such that $\la_k\geq j$, $j\in\ZZ$. The product in
 denominator is finite,
since
there are only finitely many distinct $\la'_j$ for a given choice of
$\la_k$.

\begin{lem}\label{poles}
The product $\e{n}(z_1)\e{m}(z_2)$ has poles at the points
\bse\label{pol-zer}
\bn
z=\frac{z_1}{z_2}=
\left\{\begin{array}{ll}
q^{2},\ q^{4},\ \ldots\  ,\ q^{2m},\qquad n\geq m\\
q^{2m-2n+2},\ q^{2m-2n+4},\ \ldots\  ,\ q^{2m},\qquad n<m,
\end{array}\right.
\label{polee}
\ed
and zeroes at the points
\bn\label{zeroe}
z=
\left\{\begin{array}{ll}
q^{-2n+2},\ q^{-2n+4},\ \ldots\ ,\ q^{-2n+2m},\qquad n\geq m\\
q^{-2n+2},\ q^{-2n+4},\ \ldots ,\ 1,\qquad n<m.
\end{array}\right.
\ed
Analogously, the product
$\f{n}(z_1)\f{m}(z_2)$
has
zeroes at the points
\bn \label{zerof}
z=
\left\{\begin{array}{ll}
1,\ q^{2},\ \ldots\ ,\ q^{2m-2)},\qquad n\geq m\\
q^{2m-2n},\ q^{2m-2n+2},\ \ldots ,\ q^{2m-2},\qquad n<m,
\end{array}\right.
\ed
and simple poles  at
\bn
z=
\left\{\begin{array}{ll}
q^{-2n},\ q^{-2n+2},\ \ldots\ ,\ q^{-2n+2m-2},\qquad n\geq m\\
q^{-2n},\ q^{-2n+2},\ \ldots\ ,\ q^{-2},\qquad n<m.
\end{array}\right.
\label{polef}
\ed
\ese
\end{lem}
{\it Proof.} Let us prove the part of the lemma concerning the composed
current $\e{n}(z)$. The rest can be proved analogously.
 {}From the fact that the product $e(z_1)e(z_2)$ has simple pole at the
 point
$z_1=q^2z_2$ and simple zero at $z_1=z_2$
it is clear that the product $e(q^{2k}z_1)\e{m}(z_2)$ has one simple pole at
$z_1=q^{2(m-k)}z_2$  and one simple zero at $z_1=q^{-2k}z_2$, for
$k=0,1,\ldots,n-1$. The rest poles and zeros
cancel each others. When $n<m$ this is the structure of poles and zeros
 given
in \r{polee} and \r{zeroe}. When $n> m$ there is an additional
poles/zeros cancellation. Namely, the poles at the points
$z_1=q^{-2k}z_2$, $k=0,1,\ldots,n-m-1$ cancel with the zeros at the same
points.

Put $g'(z)=\frac{q^{-2}-z}{1-q^{-2}z}$.
Then we have for all $k$, such that ${\rm max}\,(1,m-n+1)\leq k\leq m$:
\be\label{rre1}
\Res{z_1=q^{2k}z_2}\!\!\eee{n}(z_1)\eee{m}(z_2)\ddz{1}=
\lmt{z_1\to q^{2k}z_2}(1-q^{2k}z_2/z_1)
\prod_{i=0}^{k-1}\prod_{j=0}^{n-1}g'(q^{2(j-i)}z_1/z_2)
\eee{k}(z_2)\eee{n}(z_1)\eee{m-k}(q^{2k}z_2)=\\
=
\lim\limits_{z_1\to q^{2k}z_2}
(q^2-q^{2(k-1)}z_2/z_1)
\prod\limits_{i=0}^{k-2}
\prod\limits_{j=0}^{n-1}g'(q^{2j-2i}z)
\prod\limits_{j=1}^{n-1}g'(q^{2j-2k+2}z)
\eee{k}(z_2)\eee{n}(z_1)\eee{m-k}(q^{2k}z_2)=\\
=(q-q^{-1})
\prod\limits_{j=0}^{k-1}\frac{[n+j]_q[n+j+1]_q}{[j]_q[j+1]_q}
\eee{n+k}(z_2)\eee{m-k}(q^{2k}z_2)\ ,
\ee
since $g'(q^{2k})=\frac{[k+1]_q}{[k-1]_q}$ for $k=2,3,\ldots$\ .

Analogously, for  all $k$, such that ${\rm max}\,(1,n-m+1)\leq k\leq n$
\be\label{rrf1}
\Res{z_1=q^{-2k}z_2}\!\!\!\tilde{f}^{(n)}(z_1)\tilde{f}^{(m)}(z_2)\ddz{1}=
(q^{-1}-q^{})
\prod\limits_{j=0}^{k-1}\frac{[m+j]_q[m+j+1]_q}{[j]_q[j+1]_q}
\tilde{f}^{(n-k)}(z_2)\tilde{f}^{(m+k)}(q^{-2k}z_2).
\ee

These calculations together with \r{eft} give the following commutation
relations:
\be
\label{ccmme} \e{n}(z_1)\e{m}(z_2)=\prod_{k=0}^{m-1}
\prod_{l=0}^{n-1}g'(q^{2(k-l)}z)\e{m}(z_2)
\e{n}(z_1)+\\
+\sum\limits_{k=1}^{m}
[m]_q\left[{m-1\atop k}\right]_q\left[{m-1\atop k-1}\right]_q\
\delta\left(\frac{z_1}{q^{2k}z_2}\right)\e{n+k}(z_2)\e{m-k}(z_1)
\ed
Analogously, we can derive the relations
\bn
\label{ccmmf}
\f{n}(z_1)\f{m}(z_2)=\prod_{k=0}^{m-1}\prod_{l=0}^{n-1}
g(q^{2(k-l)}z)\f{m}(z_2)
\f{n}(z_1)-\\
-\sum\limits_{k=1}^{n}
[n]_q\left[{n-1\atop k}\right]_q\left[{n-1\atop k-1}\right]_q\
\delta\left(\frac{q^{2k}z_1}{z_2}\right)\f{n-k}(z_2)\f{m+k}(z_1)
\ed
 where $g(z)=\frac{q^2-z}{1-q^2z}$.
Note first, that the coefficients before delta functions are the new fields,
and second, that some of these fields  in r.h.s. of \rf{ccmme}, \rf{ccmmf}
equal zero according to the structure of zeroes \rf{pol-zer}.

Lemma \ref{poles}   shows
 also that the only (simple) poles
 of the product $\T{n}(z_1)\T{m}(z_2)$ are $z_1=q^{2m}z_2$ and
 $z_1=q^{-2n}z_2$ and from \rf{rre1} and \rf{rrf1} we deduce that
\bn
[\T{n}(z_1),\T{m}(z_2)]=
\delta(q^{2n}z_1/z_2)\T{n+m}(z_1)-
\delta(q^{-2m}z_1/z_2)\T{n+m}(z_2).
\label{T}
\ed
 so the tensor fields $\T{n}(z)$ form a closed algebra.
In particular, from \r{T} we have the following
\begin{prop}
 The total integrals of the fields $\T{n}(z)$ around infinity
$$\I{n}=\oint \T{n}(z)\dz$$
commute between themselves:
\be\label{433a}
[\I{n},\I{m}] =0.
\ee
\end{prop}

\subsection{Calculation of the contour integrals}

We want to calculate the integral
$$
\krr{n}
=\frac{1}{n!}\oint\limits_{D_n}\T{1}(z_1)\dzz{1}\cdots \T{1}(z_n)\dzz{n}
$$
over the cycles described in the previous section.
 It was noted in the proof of Proposition \r{sldvatheorem}, that for $|q|>1$
the cycle ${D}_n$ is homotopic to a product of unit circles.

 We can also deform
 it slightly to a torus $C_i=\{|z_i|=r_i\}$ in such a way, that it has no
 intersections with diagonals.
 Let us consider it as the  multiple integral
$$
\krr{n}=\frac{1}{n!}\oint\limits_{C_n}\cdots
 \oint\limits_{C_1}\T{1}(z_1)\dzz{1}\cdots \T{1}(z_n)\dzz{n}
$$
and integrate first over $z_1$ for fixed other $z_j$. We can move the
contour
 $C_1$ to infinity crossing the poles $z_1=q^2z_j$, $j=2,\ldots,n$.
The residue at the pole
 $z_1=q^2z_2$ is equal, due to \rf{rre1}, to
\footnote{Expanding the contour
to infinity we should pick up minus residue. The minus sign in the formula
\r{m-res} results that in all calculations below we avoid the appearing
of the alternating signs.}
\be\label{m-res}
-\oint\limits_{C_2}\T{2}(z_2)\dzz{2}\T{1}(z_3)\dzz{3}\cdots \T{1}(z_n)
\dzz{n}\ .
\ee
 For the
calculation of the residues at $z_1=z_3$ we first use the commutativity
  of analytical continuations of $\T{k_i}(z_i)$, which mean, in particular,
  that they commute on the  contour which does not cross their singularities
  and then repeat the same calculation as for $z_1=q^2z_2$. As a result, we
  will have $(n-1)$ integrals with different contours, obtained by
permutations
  of each other, but due to the absence of the poles near diagonal
 (see \rf{T}), all they are equivalent.

 So, we have,
$$
\krr{n}=\frac{1}{n!}\left(\ \oint\limits_{C_1}\T{1}(z_1)\dzz{1}\cdot
\oint\limits_{C_{n}}\cdots
 \oint\limits_{C_2}\T{1}(z_2)\dzz{2}\cdots \T{1}(z_{n})\dzz{n}+\right.
$$
 $$
+\left.
 (n-1)\oint\limits_{C_{n-1}}\cdots
 \oint\limits_{C_2}\T{2}(z_2)\dzz{2}\T{1}(z_3)\dzz{3}\cdots
 \T{1}(z_n)\dzz{n}\right)
$$
or
 $$
\krr{n}=\frac{1}{n}\I{1}\krr{n-1}+\frac{n-1}{n!}
 \oint\limits_{C_{n}}\cdots
 \oint\limits_{C_2}\T{2}(z_2)\dzz{2}\T{1}(z_3)\dzz{3}\cdots
\T{1}(z_n)\dzz{n})\ .
$$
 Again, integrate the second integral over $z_2$ for fixed others.
 We see, that this integral has singularities for $z_2=q^{\pm 4} z_j$.
 It is clear, that when we move the contour to infinity, we cross
  the poles $z_2=q^4z_j$, $j=3,\ldots,n$. By the same trick we have now
 $$
\krr{n}=\frac{1}{n}\I{1}\krr{n-1}
 +\frac{n-1}{n!}\oint\limits_{C_2}\T{2}(z_2)\dzz{2}\cdot
 \oint\limits_{C_{n-1}}\cdots\oint\limits_{C_{3}}
 \T{1}(z_3)\dzz{3}\cdots \T{1}(z_n)\dzz{n})+
$$
 $$
+\frac{(n-1)(n-2)}{n!}
 \oint\limits_{C_{n}}\cdots
 \oint\limits_{C_3}\T{3}(z_3)\dzz{3}\T{1}(z_4)\dzz{4}\cdots
\T{1}(z_n)\dzz{n})=
$$
 $$
=\frac{1}{n}\left(\I{1}\krr{n-1}+\I{2}\krr{n-2}+
 \frac{(n-1)(n-2)}{n!}
 \oint\limits_{C_{n}}\cdots
 \oint\limits_{C_3}\T{3}(z_3)\dzz{3}\T{1}(z_4)\dzz{4}\cdots
 \T{1}(z_n)\dzz{n}
\right).
$$
 The inductive calculation give the following recurrence relation:
 \bn
 \krr{n}=\frac{1}{n}\left(\I{1}\krr{n-1}+\I{2}\krr{n-2}+\ldots +
\I{n}\krr{0}\right)
 \label{recc}
 \ed
 with initial condition $\krr{0}=1\ot 1$.
 For the solution of this recurrence relation let us introduce
  generating functions
 $$
\ov\kr(u)=\sum_{n\geq 0}\krr{n} u^n, \qquad I(u)=\sum_{n\geq 1}\I{n}u^n.
$$
 Then the relation \rf{recc} means the following differential equation:
 \bn
 u\frac{d}{du}\ov\kr(u)=I(u)\ov\kr(u)
 \label{dif}
 \ed
 with initial condition  $\ov\kr(0)=1\ot 1$.
Its solution for the commutative variable $I(u)$ is
$$
\ov\kr(u)=\exp\sk{\int\limits_{0}^{u}\frac{I(v)}{v}dv}\ ,
$$

This proves the following
\begin{th} \label{new-for}
For the algebra $\UqdvaD$ the factor $\R$ of the universal
$R$-matrix \r{RR} can be written in the  form
\be\label{new-int}
\R=\R(1)=\exp\sk{\sum_{n\geq 1}\frac{I^{(n)}}{n}}=
\exp\sk{\sum_{n\geq 1}\oint \frac{dz}{2\pi iz}\ \frac{\T{n}(z)}{n}}.
\ee
\end{th}

\subsection{Some applications}
Let us first apply the formula \r{new-int} to the integrable representations
of $\UqdvaD$. It is known \cite{DM}, that in integrable representation of
$\UqdvaD$ of level $k$ all the currents $e^{(n)}(z)$ and $f^{(n)}(z)$ vanish
for $n>k$. It gives immediately the following statement:
\begin{prop}
The action of the universal $R$-matrix for $\UqdvaD$ in tensor product
$V\ot W$ of highest weight representations, where one of them is integrable
of level $k$, coincides with the action of
$$\K \exp\sk{\sum\limits_{n= 1}^k\oint \frac{dz}{2\pi i z}
\ \frac{\T{n}(z)}{n}}.$$
In particular, for $k=1$, it is
$$\K \exp\sk{(q^{-1}-q)\oint \frac{dz}{2\pi i z}
\ f(z)\ot e(z)}.$$
\end{prop}
Here $\K$ is given by \r{RCW1}.

We can  deduce from \r{new-int} the formula for the universal
$R$-matrix of $U_q(\sldva)$, which is well known \cite{D}. To do this, we
suppose that we act by $\kr$ on a tensor product of highest weight
representations of zero level with zero top degree and calculate the matrix
coefficient
between tensor products of the vectors of zero degree. Then in the
expression of the $R$-matrix only the terms composed from the elements from
$U_q(\sldva)$
 $e_0$, $f_0$, $q^h$ and $q^{-h}$ survive.
 What is left, is to count the
coefficients before $e_0^n$ and $f_0^n$ in the fields $e^{(n)}$ and
$f^{(n)}$. They are given in \r{n18}, \r{n19} and \r{eft}.
 Finally, we get
$$R=q^{\frac{-h\ot h}{2}}\exp\sk{\sum\limits_{n\geq 1}
\frac{(q^{-1}-q)^{2n-1}}{n[n]_q}f^n\ot
e^n}.$$
One can check that it coincides with the usual presentation of the universal
$R$-matrix via $q^2$-exponential
function using, for instance, presentation of $q$-exponent in a form of
infinite product.

Finally, we can use factorized expression of \cite{DK} for the universal
$R$-matrix and substitute there \r{new-int} in every factor. As a result, we
get an expression of the form, which depends on a fixed normal ordering
 $>$ of the
system $\Delta_+$  of positive roots (or, equivalently, on the reduced
 decomposition of the longest element of the Weyl group) of simple Lie algebra
$\g$ of simple laced type:

$$
\kr=\K
\stackrel{\longrightarrow}{\prod_{\gamma\in\Delta_+}}
\exp\sk{\sum_{n\geq1} \oint \frac{dz}{2\pi i z}\ \frac{t_\gamma^{(n)}(z)}{n}}
$$
where $\Delta_+$ is the system of positive roots of simple Lie algebra $\g$,
the currents $e^{}_{\ga}(z)$ and $e^{}_{-\ga}(z)$ are built as in \cite{DK},
$t^{}_{\ga}(z)=(q^{-1}-q)e^{}_{-\ga}(z)\ot e^{}_{\ga}(z)$
 and the currents $t^{(n)}_{\ga}(z)$ and    are constructed
from $e^{}_{\ga}(z)$ and $e^{}_{-\ga}(z)$ according to the rules
\r{eft}. Again, we can restrict ourselves to zero level matrix coefficients
and recover precisely multiplicative formula for the $R$-matrix of
$U_q(\g)$, see, e.g., \cite{KT}.

We can also repeat our considerations for Yangians and elliptic algebras, as it
was done in \cite{DK} and in \cite{DKKP}.
Another applications are given in \cite{DKP}

\section*{Acknowledgment}

\qquad S. Khoroshkin
is grateful to B.Feigin and Yu.Nerertin for the discussions
 and thanks especially  A.Rosly for the explanations in
 the theory of residues.
He also thanks the University of
Chicago and the University of Cincinnati for hospitality. Part of the work
was done there. S. Pakuliak thanks the Physikalisches Institut,
Universit\"at Bonn;
CERN Theory Division and  Universite d'Angers for hospitality.

The work of S.Kh. and of S.P.
 was supported in part by the grants INTAS OPEN 97-01312,
RFBR-CNRS grant PICS N 608/RFBR 98-01-22033.
S.Kh. was supported also  by RFBR grant 98-01-00303 and S.P.
by grants RFBR 97-01--1041  and grant of Heisenberg-Landau
program HL-99-12.

\end{document}